%% file: skm.tex
\documentclass[12pt]{amsart}
\renewcommand{\bar}{\overline}

\newcommand{\pa}{\partial}
\newcommand{\ph}{\varphi}

\usepackage{amsmath,amsthm,amscd}

\title
[]{A Note on  Special K\"ahler Manifolds}
\author[]{Zhiqin Lu}
\date{\today}
\address[Zhiqin Lu]
{Department of Mathematics\\
Columbia University\\
NY, NY 10027}
\email{lu@cpw.math.columbia.edu}

\newtheorem{theorem}{Theorem}

\newtheorem{prop}{Proposition}

\newtheorem{definition}{Definition}

\theoremstyle{remark}

\begin{document}
\maketitle

\input{skm10.tex}

\bibliographystyle{plain}
\bibliography{bib}

\end{document}

%% file: skm10.tex
The base space of an algebraically completely integrable Hamiltonian
system acquires a rather special differential-geometric structure which
plays an important role in modern physical theories such as Seiberg-Witten
theory. This structure was formalised by D. Freed~\cite{DSF} as a 
\textsl{Special K\"ahler manifold}. In that paper he conjectured that
there are no compact special K\"ahler manifolds other than flat ones. 
In this paper, we prove that there are no nonflat
complete special K\"ahler manifolds, thus verifying the conjecture.

\begin{definition}
Let $M$ be a K\"ahler manifold with K\"ahler form $\omega$. A special
K\"ahler structure on $M$ is a real flat torsionfree symplectic
connection $\nabla$ satisfying
\[
d_\nabla I=0
\]
where $I\in\Omega^1(M,TM)$ is the complex structure on $M$.
\end{definition}

The following property is true for any special K\"ahler manifold.

\begin{theorem}\label{nfud}
Any special K\"ahler manifold has nonnegative Ricci curvature. Moreover,
if the scalar curvature is identically zero, then the curvature tensor
itself
vanishes.
\end{theorem}

{\bf Proof:} We take $(z^1,\cdots,z^n)$ to be a special
coordinate system. Under this coordinate system, the K\"ahler form
can be represented as
\[
\omega=\frac{\sqrt{-1}}{2} h_{i\bar j} dz^i\wedge d\bar z^j
=\frac{\sqrt{-1}}{2} Im\, u_{ij} dz^i\wedge d\bar z^j
\]
where $u$ is a local holomorphic function 
~\cite[page 5]{DSF} and $u_{ij}$ is defined as 
$\frac{\pa^2 u}{\pa z^i\pa z^j}$. Moreover, we define $u_{i_1\cdots i_n}$ to be
$\frac{\pa^n u}{\pa z^1\cdots\pa z^n}$.

It is then a straightforward computation that
\[
R_{i\bar jk\bar l}=-\frac 14 h^{m\bar n} u_{ikm}\bar {u_{jln}}
\]
where our definition for the curvature tensor is
\[
R_{i\bar jk\bar l}=\frac{\pa^2 h_{i\bar j}}{\pa z^k\pa\bar z^l}
-h^{m\bar n}\frac{\pa h_{m\bar j}}{\pa \bar z^l}\frac{\pa h_{i\bar n}}{\pa z^k}
\]

Let $F$ be the cubic form defined as
\[
F=-\omega(\pi^{(1,0)}, \nabla \pi^{(1,0)})
\]
where $\pi^{(1,0)}\in\Omega^{1,0}(T_CM)$ is the  projection onto the $(1,0)$
part of the complexified tangent bundle.
$F$ is a global section of the bundle $Sym^3 T^*M$.
Locally $F_{ijk}=u_{ijk}$. Thus we
have
\begin{equation}\label{xx}
R_{i\bar jk\bar l}=-\frac 14 h^{m\bar n} F_{ikm}
\bar{F_{jln}}
\end{equation}

We remark that the  above identity is valid for any local holomorphic coordinate
system,
because both sides are tensors.
From~\eqref{xx} we see that the Ricci curvature is nonnegative.
Furthermore, if 
the scalar curvature $\rho$ 
\begin{equation}\label{1}
\rho=\frac 14 h^{i\bar i_1} h^{j\bar j_1}h^{k\bar k_1}
F_{ijk}\bar{F_{i_1j_1k_1}}
\end{equation}
is identically zero, then the curvature tensor vanishes.

\qed

\begin{theorem}
If $(M,\omega)$ is a complete special K\"ahler manifold, then the
curvature  tensor
vanishes.
\end{theorem}

{\bf Proof:}
Define
\begin{equation}\label{2}
F_{ijk,l}=
\pa_l F_{ijk} -\Gamma^m_{il} F_{mjk}
-\Gamma^m_{jl} F_{imk}
-\Gamma^m_{kl} F_{ijm}
\end{equation}
to be the covariant derivative of $F_{ijk}$. By the Bochner formula
\begin{align*}
&\Delta\rho=\frac 14
 h^{i\bar i_1} h^{j\bar j_1}h^{k\bar k_1}h^{l\bar l_1}
F_{ijk,l}\bar{F_{i_1j_1k_1,l_1}}\\
&+\frac 34 h^{m_1\bar i_1} h^{i\bar m}
h^{j\bar j_1} h^{k\bar k_1} R_{m_1\bar m} F_{ijk}\bar{F_{i_1j_1k_1}}
\end{align*}
where $R_{m_1\bar m}=-h^{\alpha\bar\beta} R_{\alpha\bar\beta m_1\bar m}$
is the Ricci tensor.
In particular, we have
\begin{equation}\label{yy}
\Delta\rho\geq  3h^{m_1\bar i_1} h^{i\bar m}R_{i\bar i_1} R_{m_1\bar m}
\geq \frac 3n\rho^2
\end{equation}

We use the  version of the generalized maximal principle
in~\cite[page 582, Lemma 1.1]{TY1}. 
The original statement in~\cite{TY1} is quite general. We rewrite it in
the following (simplified) way:

\begin{prop}
Let $M$ be a complete K\"ahler manifold of nonnegative Ricci curvature.
Let $\ph$ be a nonnegative function satisfying
\[
\Delta\ph\geq C_1\ph^\alpha-C_2\ph-C_3
\]
where $\alpha>1$, $C_1>0$, $C_2,C_3\geq 0$ are constants. Then
\[
sup\,\ph\leq Max\{1,\left(\frac{C_2+C_3}{C_1}\right)^{1/\alpha}\}
\]
\end{prop}

\qed

For any positive number $a$, from ~\eqref{yy} we see that
$\Delta (a\rho)\geq \frac 3{na} (a\rho)^2$. By the above proposition we
see
that $a\rho\leq 1$. Letting $a\rightarrow \infty$ we see that $\rho\equiv
0$. The 
theorem thus follows from
Theorem~\ref{nfud}.
 \qed

{\bf Acknowledgment.} The author thanks Professor G. Tian for the 
support during the preparation of this paper.
He  thanks Professor D. Freed for some very useful conversations
about
this topic, from which he learned a lot about the special K\"ahler
geometry. He also thanks 
Professor N. Hitchin for the interest of the paper and suggestions
which make the paper in its current form.